\documentclass[10pt,technote]{IEEEtran}
\usepackage{graphicx}          
\usepackage{amsmath}           
\usepackage{amssymb}
\usepackage{cite}
\usepackage{algorithm}
\usepackage{enumerate}
\usepackage{algpseudocode}
\usepackage[latin1]{inputenc}  

\usepackage{amsthm}

\theoremstyle{definition}

\theoremstyle{remark}
\theoremstyle{plain}

\theoremstyle{remark} \newtheorem{note}{Remark}

\DeclareMathOperator{\E}{E}

\newcommand{\mbf}[1]{\mathbf{#1}}
\newcommand{\mbs}[1]{\boldsymbol{#1}}

\begin{document}

\title{Cramér-Rao bound analog of Bayes rule}
\author{Dave Zachariah and Petre Stoica}

\maketitle



Estimation of multiple parameters is a common task in signal
processing. The Cramér-Rao bound (CRB) sets a statistical lower
limit on the resulting errors when estimating parameters from
a set of random observations. It can be understood as a fundamental
measure of parameter uncertainty \cite{Cramer1946_contribution,vanTrees2013_detection}. As a
general example, suppose
$\mbs{\theta}$ denotes the vector of sought parameters and that the random
observation model can be written as
\begin{equation}
\mbf{y} = \mbf{x}_\theta + \mbf{w},
\label{eq:additivemodel}
\end{equation}
where $\mbf{x}_{\theta}$ is a function or signal parameterized by
$\mbs{\theta}$ and $\mbf{w}$ is a zero-mean Gaussian noise
vector. Then the CRB for $\mbs{\theta}$ has the following notable properties:
\begin{enumerate}[i)]
\item For a fixed $\mbs{\theta}$, the CRB for
$\mbs{\theta}$ decreases as the dimension of $\mbf{y}$ increases.

\item For a fixed $\mbf{y}$, if additional parameters
  $\tilde{\mbs{\theta}}$ are estimated then the CRB for $\mbs{\theta}$ increases as the dimension of
$\tilde{\mbs{\theta}}$ increases.

\item  If adding a set of
observations $\tilde{\mbf{y}}$ requires estimating additional
parameters $\tilde{\mbs{\theta}}$, then the CRB for
$\mbs{\theta}$ decreases as the dimension of  $\tilde{\mbf{y}}$ increases,
provided the dimension of $\tilde{\mbs{\theta}}$ does not exceed that
of $\tilde{\mbf{y}}$
\cite{Stoica&Li1996_study}. This property implies both i) and ii) above.

\item Among all possible distributions of $\mbf{w}$ with a fixed covariance
  matrix, the CRB for $\mbs{\theta}$ attains its maximum when $\mbf{w}$
  is Gaussian, i.e., the Gaussian scenario is the `worst-case' for estimating $\mbs{\theta}$ \cite{Stoica&Babu2011_gaussian,ParkEtAl2013_gaussian,SteinEtAl2014_lower}.
\end{enumerate}


In this lecture note, we show a general property of the CRB
that quantifies the interdependencies between the parameters in
$\mbs{\theta}$. The presented result is valid for more general models than
\eqref{eq:additivemodel} and also generalizes the result in
\cite{DAmico2014_reciprocity} to vector parameters. It will be illustrated via two examples.


\section{Relevance}
In probability theory, the chain rule and Bayes rule are useful
tools to analyze the statistical interdependence between multiple
random variables and to derive tractable expressions for their
distributions. In this lecture note, we provide analogs of the chain
rule and Bayes rule for the Cramér-Rao bound associated with multiple
parameters. The results are particularly useful when estimating parameters of interest in the presence
of nuisance parameters.

\section{Prerequisities}

The reader needs basic knowledge about linear algebra,
elementary probability theory, and statistical signal processing.

\section{Preliminaries}
We will consider a general scenario in which we observe an $n
\times 1$ random vector $\mbf{y}$. Its probability density function (pdf) $p(\mbf{y};
\mbs{\theta})$ is
parameterized by a $k \times 1$ deterministic vector $
\mbs{\theta}$. The goal is to estimate $
\mbs{\theta}$, or subvectors
of $\mbs{\theta}$, given $\mbf{y}$. 

Let $l(\mbs{\theta}) \triangleq \ln p(\mbf{y};\mbs{\theta})$ denote the log-likelihood
function and let $\hat{\mbs{\theta}}$ be any unbiased estimator. Then the mean square error (MSE) matrix $\mbf{P}_{\hat{\theta}} \triangleq \E[ (\mbs{\theta}
- \hat{\mbs{\theta}})(\mbs{\theta}
- \hat{\mbs{\theta}})^* ]$ is bounded from below by the inverse of the Fisher
information matrix $\mbf{J}_\theta \triangleq -\E[\partial^2_\theta l(\mbs{\theta}) ]$, where $\partial^2_\theta$
denotes the second-order differential or Laplacian operator with
respect to $\mbs{\theta}$. That is, $\mbf{P}_{\hat{\theta}}
\succeq \mbf{J}^{-1}_\theta$, assuming from hereon that $\mbf{J}_\theta$ is
nonsingular. This is the Cramér-Rao inequality \cite{vanTrees2013_detection,Soderstrom&Stoica1988_system,Kay1993_fundamentals}.

The determinant of the MSE-matrix, $|\mbf{P}_{\hat{\theta}}|$, is a
scalar measure of the error magnitude. For unbiased estimators,
$|\mbf{P}_{\hat{\theta}}|$ equals the `generalized variance' of errors \cite{Anderson2003_introduction}. By defining $\text{CRB}(\mbs{\theta}) \triangleq
|\mbf{J}^{-1}_\theta|$, the generalized error variance is bounded by
\begin{equation*}
|\mbf{P}_{\hat{\theta}}| \geq \text{CRB}(\mbs{\theta}).
\end{equation*}

In the following we are interested in subvectors or elements of
$\mbs{\theta}$. Letting $\mbs{\theta} = [\mbs{\alpha}^\top \;
\mbs{\beta}^\top]^\top$, we can write the Fisher information
matrix in block-form,
\begin{equation}
\mbf{J}_\theta = -\E 
\begin{bmatrix}
\partial^2_\alpha l( \mbs{\alpha}, \mbs{\beta}) & \partial_\alpha \partial_\beta l(
\mbs{\alpha}, \mbs{\beta}) \\
\partial_\beta \partial_\alpha l(
\mbs{\alpha}, \mbs{\beta}) & \partial^2_\beta l(
\mbs{\alpha}, \mbs{\beta}) 
\end{bmatrix}  
= \begin{bmatrix}
\mbf{J}_\alpha & \mbf{J}_{\alpha \beta} \\
\mbf{J}_{\beta \alpha} & \mbf{J}_\beta
\end{bmatrix}.
\label{eq:jointFIM}
\end{equation}


\section{Main result}
Let $\mbf{a}$ and $\mbf{b}$ be two random vectors. Two useful rules
in probability theory are the chain rule
\begin{equation}
p(\mbf{a},\mbf{b}) = p(\mbf{a}|\mbf{b})p(\mbf{b})
\label{eq:chainrule_prob}
\end{equation}
and Bayes rule
\begin{equation}
p(\mbf{a}) = \frac{p(\mbf{b})}{p(\mbf{b} | \mbf{a})} p(\mbf{a}| \mbf{b} ).
\label{eq:bayesrule_prob}
\end{equation}

Now consider two parameter vectors $\mbs{\alpha}$ and $\mbs{\beta}$. When both are
unknown, their
\emph{joint} Cramér-Rao bound is given by
\begin{equation}
\text{CRB}(\mbs{\alpha}, \mbs{\beta}) = \left |
\begin{bmatrix}
\mbf{J}_\alpha & \mbf{J}_{\alpha \beta} \\
\mbf{J}_{\beta \alpha} & \mbf{J}_\beta
\end{bmatrix}^{-1}
\right |.
\label{eq:crb_joint}
\end{equation}
The bound for $\mbs{\alpha}$ with \emph{known} $\mbs{\beta}$ is simply
\begin{equation}
\text{CRB}(\mbs{\alpha} | \mbs{\beta}) = \left | \mbf{J}^{-1}_{\alpha} \right |,
\label{eq:crb_alphabeta}
\end{equation}
and the bound for $\mbs{\alpha}$ with \emph{unknown} $\mbs{\beta}$ is
\begin{equation}
\text{CRB}(\mbs{\alpha}) =  \left | (\mbf{J}_{\alpha} - \mbf{J}_{\alpha \beta}
  \mbf{J}^{-1}_{\beta} \mbf{J}_{\beta \alpha} )^{-1} \right |.
\label{eq:crb_alpha}
\end{equation}
(\eqref{eq:crb_alpha} follows by evaluating the inverse in \eqref{eq:crb_joint} and
extracting the upper-left block corresponding to
$\mbs{\alpha}$.) Eqs. \eqref{eq:crb_alphabeta} and \eqref{eq:crb_alpha}
are the respective CRB analogs of conditional and marginal distributions for random
variables.

By applying the Schur determinant formula \cite{Horn&Johnson1990_matrix,Soderstrom&Stoica1988_system}
\begin{equation*}
\begin{split}
\left |
\begin{bmatrix}
\mbf{J}_\alpha & \mbf{J}_{\alpha \beta} \\
\mbf{J}_{\beta \alpha} &\mbf{J}_\beta
\end{bmatrix} \right | &= |\mbf{J}_\alpha||\mbf{J}_{\beta} - \mbf{J}_{\beta \alpha}
\mbf{J}^{-1}_{\alpha} \mbf{J}_{\alpha \beta} | \\
&=  |\mbf{J}_\beta||\mbf{J}_{\alpha} - \mbf{J}_{\alpha \beta}
\mbf{J}^{-1}_{\beta} \mbf{J}_{\beta \alpha} |,
\end{split}
\end{equation*}
along with $|\mbf{J}^{-1}| = |\mbf{J}|^{-1}$, to
\eqref{eq:crb_joint}, \eqref{eq:crb_alphabeta} and \eqref{eq:crb_alpha}, we
can now state the Cramér-Rao bound analogs of the chain rule \eqref{eq:chainrule_prob},
\begin{equation}
\begin{split}
\text{CRB}(\mbs{\alpha}, \mbs{\beta}) = \text{CRB}(\mbs{\alpha} | \mbs{\beta}) \text{CRB}(\mbs{\beta})
\end{split}
\label{eq:chainrule}
\end{equation}
and of Bayes rule \eqref{eq:bayesrule_prob},
\begin{equation}
\begin{split}
\boxed{\text{CRB}(\mbs{\alpha}) = \frac{\text{CRB}(\mbs{\beta})}{\text{CRB}(\mbs{\beta} |
  \mbs{\alpha})} \text{CRB}(\mbs{\alpha} | \mbs{\beta})}.
\end{split}
\label{eq:bayesrule}
\end{equation}
The results are of course symmetric, i.e., one can interchange
$\mbs{\alpha}$ and $\mbs{\beta}$.

From \eqref{eq:chainrule} we see that the \emph{joint} error bound for $\mbs{\alpha}$ and $\mbs{\beta}$ equals the error
bound for $\mbs{\alpha}$, when $\mbs{\beta}$ is known, multiplied by the error
bound for $\mbs{\beta}$. More interestingly, \eqref{eq:bayesrule} tells us
that the error bound for $\mbs{\alpha}$ is equal to the
bound for $\mbs{\alpha}$ when $\mbs{\beta}$ is known, multiplied by a factor, viz.
$\text{CRB}(\mbs{\beta})/\text{CRB}(\mbs{\beta} | \mbs{\alpha}) \geq 1$, that quantifies
the influence of $\mbs{\beta}$ on one's ability to estimate
$\mbs{\alpha}$. 

\begin{note}
The rules can be applied to cases with any number of additional parameters, besides $\mbs{\alpha}$ and $\mbs{\beta}$. Consider
for instance the case of $\mbs{\alpha}$, $\mbs{\beta}$ and $\mbs{\gamma}$, where $\mbs{\gamma}$ is an
unknown nuisance parameter. Then applying the chain rule twice yields
\begin{equation}
\begin{split}
\text{CRB}(\mbs{\alpha}, \mbs{\beta}, \mbs{\gamma}) &= \text{CRB}(\mbs{\gamma}| \mbs{\alpha}, \mbs{\beta})
\text{CRB}(\mbs{\alpha} | \mbs{\beta}) \text{CRB}(\mbs{\beta})\\
&= \text{CRB}(\mbs{\gamma}| \mbs{\alpha}, \mbs{\beta}) \text{CRB}(\mbs{\beta} | \mbs{\alpha}) \text{CRB}(\mbs{\alpha})
\end{split}
\label{eq:chainrule2}
\end{equation}
where the factors without $\mbs{\gamma}$ signify that the nuisance parameter
is \emph{unknown}. Combining the two expressions in
\eqref{eq:chainrule2} yields the analog of Bayes rule
\eqref{eq:bayesrule} for any number of additional
parameters.

The joint error bound for a set of parameters $\mbs{\alpha}_1, \mbs{\alpha}_{2}, \mbs{\alpha}_{3},
\dots$ can be similarly decomposed by a recursive application
of the chain rule in order to analyze their interdependency and its
impact on estimation.
\end{note}

\begin{note}
The CRB analog of Bayes rule \eqref{eq:bayesrule} generalizes
the result in \cite{DAmico2014_reciprocity} which concerns only
scalar parameters $\alpha$ and $\beta$ amid a vector of nuisance parameters $\mbs{\gamma}$. Our proof of \eqref{eq:bayesrule} is also more direct than in \cite{DAmico2014_reciprocity}.
\end{note}

\begin{note}
These results are also applicable to the posterior,
or Bayesian, Cramér-Rao bound (PCRB), in which $\mbs{\theta}$ is modeled as a random
variable with a prior distribution. The PCRB is valid for the entire class of estimators $\hat{\mbs{\theta}}$,
whether biased or not \cite{vanTrees2013_detection}. The posterior Cramér-Rao inequality is then $\mbf{P}_{\hat{\theta}} \succeq \mbf{J}^{-1}_{\theta}$, where $\mbf{J}_{\theta} \triangleq -\E[ \partial^2_\theta \ln p(\mbf{y}, \mbs{\theta}) ]$ is the Bayesian Fisher information matrix, $p(\mbf{y},\mbs{\theta})$ is the joint pdf and the expectation is with respect to this pdf.  Letting $\mbs{\theta} = [\mbs{\alpha}^\top \;
\mbs{\beta}^\top]^\top$, the matrix can be partitioned correspondingly,
\begin{equation*}
\mbf{J}_{\theta} = \begin{bmatrix} \mbf{J}_\alpha & \mbf{J}_{\alpha \beta} \\ \mbf{J}_{\beta \alpha} & \mbf{J}_{ \beta} \end{bmatrix},
\end{equation*}
and thereby the results \eqref{eq:chainrule}, \eqref{eq:bayesrule} and \eqref{eq:chainrule2} can be applied to the PCRB as well.
\end{note}

\section{Examples}

Next, we illustrate via two examples how a decomposition like \eqref{eq:bayesrule} can be used for analysis. The examples show that, by quantifying the
impact of nuisance parameters, it is possible to study the trade-off between the gain of obtaining them through independent side information versus estimating them jointly with the parameters of interest.

\subsection{Linear mixed model}

Consider a linear model 
\begin{equation*}
\mbf{y}= \mbf{A} \mbf{x} + \mbf{B} \mbf{z} + \mbf{w} \in \mathbb{R}^n,
\end{equation*}
where $\mbf{w}$ is Gaussian noise with covariance matrix
$v\mbf{I}$, and $\mbf{x} \in \mathbb{R}^{k_x}$ and $\mbf{z}\in \mathbb{R}^{k_z}$ are
unknown parameters. The matrices are known and $\text{rank}([ \mbf{A}
\;  \mbf{B}]) = k_x + k_z < n$, which implies that the parameters
$\mbf{x}$ and $\mbf{z}$ are embedded into two distinct range spaces,
$\mathcal{R}(\mbf{A})$ and $\mathcal{R}(\mbf{B})$,
respectively. Here $\mathcal{R}(\mbf{A})$ denotes the linear subspace
spanned by the columns of $\mbf{A}$. Under these conditions the
joint Fisher information matrix equals \cite{Kay1993_fundamentals}
\begin{equation*}
\begin{bmatrix}
\mbf{J}_x & \mbf{J}_{x z} & \mbf{J}_{x v} \\
\mbf{J}_{z x} &\mbf{J}_z & \mbf{J}_{z v} \\
\mbf{J}_{v x} &\mbf{J}_{v z} & {J}_{v} 
\end{bmatrix} = \frac{1}{v}
\begin{bmatrix}
\mbf{A}^\top \mbf{A} & \mbf{A}^\top \mbf{B}& \mbf{0} \\
\mbf{B}^\top \mbf{A} &\mbf{B}^\top \mbf{B} & \mbf{0} \\
\mbf{0} &\mbf{0} & \frac{n}{2v} 
\end{bmatrix}.
\end{equation*}
From this expression, we see that the bound for $v$ is independent 
of that for $\mbf{x}$ and $\mbf{z}$. That is, $\text{CRB}(\mbf{x},\mbf{z},v) =
\text{CRB}(\mbf{x},\mbf{z})\text{CRB}(v)$. This is a CRB analog of the
independence for random variables. Furthermore, we obtain $\text{CRB}(\mbf{z} |
  \mbf{x}) = |\mbf{J}^{-1}_z| = |v (\mbf{B}^\top \mbf{B})^{-1} | = v^{k_z}|
  \mbf{B}^\top \mbf{B}|^{-1} $ and
$\text{CRB}(\mbf{z}) = |(\mbf{J}_z - \mbf{J}_{zx}\mbf{J}^{-1}_x \mbf{J}_{xz})^{-1} | = |v (\mbf{B}^\top \mbf{B} - \mbf{B}^\top \mbf{A}
(\mbf{A}^\top \mbf{A})^{-1} \mbf{A}^\top \mbf{B} )^{-1} | = v^{k_z}|
  \mbf{B}^\top \mbs{\Pi}^\perp_{\mbf{A}} \mbf{B}|^{-1} $, where
  $\mbs{\Pi}^\perp_{\mbf{A}}$ is the projector onto the orthogonal
   complement of $\mathcal{R}(\mbf{A})$.

The increase in the error bound for $\mbf{x}$ due to the lack of information about $\mbf{z}$ can now be quantified using \eqref{eq:bayesrule}
\begin{equation}
\begin{split}
\text{CRB}(\mbf{x}) 
&= \frac{|
  \mbf{B}^\top \mbf{B}|}{|\mbf{B}^\top \mbs{\Pi}^\perp_{\mbf{A}}
  \mbf{B}|} \text{CRB}(\mbf{x} | \mbf{z}),
\end{split}
\label{eq:crb_linearmodel}
\end{equation}
where the factor $|\mbf{B}^\top \mbs{\Pi}^\perp_{\mbf{A}}
\mbf{B}|$ measures the alignment of $\mathcal{R}(\mbf{A})$ and
$\mathcal{R}(\mbf{B})$. When the range spaces are orthogonal we have that
$|\mbf{B}^\top \mbs{\Pi}^\perp_{\mbf{A}}   \mbf{B}| = |\mbf{B}^\top
\mbf{B} |$, and by \eqref{eq:crb_linearmodel} the bound for $\mbf{x}$
is unaffected by one's ignorance about $\mbf{z}$. In scenarios where it is
possible to obtain $\mbf{z}$ through additional side-information or calibration, instead
of estimation, the cost can be weighed against the reduction of the error bound for $\mbf{x}$ by the given factor $|\mbf{B}^\top \mbs{\Pi}^\perp_{\mbf{A}}
  \mbf{B}| / |
  \mbf{B}^\top \mbf{B}|$. 
  
This example has illustrated the interdependencies between the unknown parameters
$\mbf{x}$, $\mbf{z}$ and $v$. Next we consider an example where 
the unknown parameters become
asympotically independent as the number of samples $n$ grows large.

\subsection{Sine-wave fitting}

Sine-wave fitting is a problem that arises in system testing, for example of waveform recorders, and the IEEE
Standard 1057 formalizes procedures to do so (\cite{Andersson&Handel2006_ieee} and references therein).

Consider $n$ uniform samples of a sinusoid in noise
\begin{equation*}
\begin{split}
y(k) &= \alpha \sin( \omega k + \phi) + C + w(k), \\
\end{split}
\end{equation*}
where $w(k)$ is a Gaussian white noise process with variance $v$ and $k=0,\dots,n-1$.  The amplitude
$\alpha$ and phase $\phi$ of the sinusoidal signal, along with the
offset $C$, are of interest. In certain cases
the frequency $\omega$ of the test signal may be obtained separately from the estimation of $\alpha, \phi$ and $C$. For simplicity, we first consider an alternative
parameterization of the sinusoid, namely:
$\alpha \sin(
\omega k + \phi)  = A \cos( \omega k) + B \sin(\omega k) $,
where $A = \alpha \sin(\phi)$ and $B = \alpha \cos(\phi)$. The
parameters are $\mbs{\theta} = [A \: B \: C
\: \omega \: v]^\top$.

As shown in \cite{Andersson&Handel2006_ieee},
the Fisher information
matrix can be decomposed into
$\mbf{J}_\theta = \bar{\mbf{J}}_\theta +
\tilde{\mbf{J}}_\theta$, where 
\begin{equation*}
\begin{split}
 \bar{\mbf{J}}_\theta &= 
\begin{bmatrix}
\bar{J}_{A} & \bar{J}_{AB} & \bar{J}_{AC} & \bar{J}_{A\omega} & \bar{J}_{A v} \\
\bar{J}_{BA} & \bar{J}_{B} & \bar{J}_{BC} & \bar{J}_{B\omega} & \bar{J}_{B v} \\
\bar{J}_{CA} & \bar{J}_{CB} & \bar{J}_{C} & \bar{J}_{C\omega} & \bar{J}_{C v} \\
\bar{J}_{\omega A} & \bar{J}_{\omega B} & \bar{J}_{\omega C} & \bar{J}_{\omega} & \bar{J}_{\omega  v} \\
\bar{J}_{vA} & \bar{J}_{v B} & \bar{J}_{v C} & \bar{J}_{v \omega} & \bar{J}_{v}
\end{bmatrix}
\\
&=\frac{1}{2v}
\begin{bmatrix}
n & 0 & 0   & - \frac{Bn^2}{2} & 0 \\
0 & n & 0   & \frac{An^2}{2} & 0 \\
0 & 0 & 2n & 0                    & 0 \\
\frac{-Bn^2}{2} & \frac{An^2}{2} & 0 & \frac{(A^2+B^2)n^3}{3} & 0 \\
0 & 0 & 0 & 0 & \frac{n}{v} 
\end{bmatrix}
\end{split}
\end{equation*}
contains the dominant terms and
$\tilde{\mbf{J}}_\theta$ contains the remainder, so that
$\mbf{J}^{-1}_\theta \simeq \bar{\mbf{J}}^{-1}_\theta $ for large
$n$. Using this approximation we now analyze the bounds for $A$, $B$
and $C$ by application of \eqref{eq:bayesrule}.

First, let $\mbs{\theta}' = [A \; B \; C \; v]^\top$ be the
parameter vector without $\omega$. Then
\begin{equation*}
\begin{split}
\text{CRB}(\omega ) &= |J_\omega - \mbf{J}_{\omega \theta'}
\mbf{J}^{-1}_{\theta'} \mbf{J}_{\theta' \omega} |^{-1} \\
&\simeq 2v\left( \frac{(A^2 + B^2)n^3}{3} - \frac{(A^2+B^2)n^3}{4}
\right)^{-1} \\
&= \frac{2v}{n^3}\frac{12}{(A^2 + B^2)}.
\end{split}
\end{equation*}
Second, let $\mbs{\theta}'' = [ B \; C \; v]^\top$ be the
parameter vector without $\omega$ and $A$. Then
\begin{equation*}
\begin{split}
\text{CRB}(\omega | A ) &= |J_\omega - \mbf{J}_{\omega \theta''}
\mbf{J}^{-1}_{\theta''} \mbf{J}_{\theta'' \omega} |^{-1} \\
&\simeq 2v\left( \frac{(A^2 + B^2)n^3}{3} - \frac{A^2n^3}{4}
\right)^{-1} \\
&= \frac{2v}{n^3}\frac{12}{(A^2 + B^2) + 3B^2}.
\end{split}
\end{equation*}
Thus $\text{CRB}(\omega ) / \text{CRB}(\omega | A ) = 1 + 3B^2/(A^2 +
B^2) \in [1,4]$. Note that $\mbf{J}_{\theta'}$ and
$\mbf{J}_{\theta''}$ are diagonal, making their inverses
particularly easy to compute. Applying \eqref{eq:bayesrule} we obtain
\begin{equation*}
\begin{split}
\text{CRB}(A) &\simeq \left(1 + \frac{3B^2}{A^2 + B^2} \right) \text{CRB}(A|\omega) \\
\text{CRB}(B) &\simeq \left(1 + \frac{3A^2}{A^2 + B^2} \right)  \text{CRB}(B|\omega) \\
\text{CRB}(C) &\simeq \text{CRB}(C|\omega),
\end{split}
\end{equation*}
where the bounds for $B$ and $C$ are derived in a similar manner as for $A$. This shows that the bound for the offset $C$ becomes independent of the
knowledge of the frequency $\omega$ as $n$ increases, while the bounds
for $A$ and $B$ are inflated by factors ranging between 1 and 4 due to one's
ignorance about $\omega$.

When considering the original parameterization $\mbs{\vartheta} = [\alpha \;
\phi \; C \; \omega \; v]^\top$ there
exists an invertible relation, $\mbs{\vartheta} =
\mbf{g}(\mbs{\theta}) = [ \sqrt{A^2 + B^2} \:
\text{arctan}(\frac{A}{B}) \: C \: \omega \: v]^\top$. Therefore we have that
$\mbf{J}^{-1}_\vartheta = \partial_\theta \mbf{g}(\mbs{\theta})
\mbf{J}^{-1}_\theta \partial_\theta
\mbf{g}(\mbs{\theta})^\top$ \cite{vanTrees2013_detection}, where $\partial_\theta$ denotes the
first-order differential or gradient with respect to $\mbs{\theta}$ and
\begin{equation*}
\partial_\theta \mbf{g}(\mbs{\theta})=
\begin{bmatrix}
\sin\phi & \cos \phi & \mbf{0} \\
\frac{\cos \phi}{\alpha} & -\frac{\sin \phi}{\alpha} & \mbf{0} \\
\mbf{0} & \mbf{0} & \mbf{I}
\end{bmatrix}.
\end{equation*}
Exploiting the approximation $\mbf{J}^{-1}_\theta \simeq \bar{\mbf{J}}^{-1}_\theta$ once again, one obtains \cite{Andersson&Handel2006_ieee} 
\begin{equation*}
\begin{split}
\text{CRB}(\alpha) &\simeq  \text{CRB}(\alpha|\omega) \\
\text{CRB}(\phi) &\simeq 4 \text{CRB}(\phi|\omega).
\end{split}
\end{equation*}
This shows that in large samples the error bound for the
amplitude $\alpha$ also becomes independent of knowledge about the frequency $\omega$, whereas 
not knowing $\omega$ inflates the bound for the
phase $\phi$ by a factor of 4.

For large data records, the cost of pre-calibrating the frequency can
be weighed against a reduction of the error bound for the phase, while the
error bounds for the amplitude and offset will not be improved.

\section{What we have learned}

An analog of Bayes rule for the Cramér-Rao bound has been
derived. This analogous rule enables a formalized decomposition and quantification of the
mutual dependencies between multiple unknown parameters. The use of
the rule was
illustrated in two estimation problems.

\section{Authors}
Dave Zachariah (dave.zachariah@it.uu.se) is a researcher and Petre Stoica
(ps@it.uu.se) is a professor, both are at the Department of Information 
Technology, Uppsala University, Uppsala, 
Sweden. 

\bibliographystyle{ieeetr}
\bibliography{refs_crbrules}

\end{document}